# Multivariate Regular Variation on Cones: Application to Extreme Values, Hidden Regular Variation and Conditioned Limit Laws

Sidney I. Resnick,
School of Operations Research and Information Engineering
Cornell University
Ithaca, NY 14853

We attempt to bring some modest unity to three subareas of heavy tail analysis and extreme value theory:

- limit laws for componentwise maxima of iid random variables;
- hidden regular variation and asymptotic independence;
- conditioned limit laws when one component of a random vector is extreme.

The common theme is multivariate regular variation on a cone and the three cases cited come from specifying the cones $[0,\infty]^d \setminus \{\mathbf{0}\}$; $(0,\infty]^d$; and $[0,\infty] \times (0,\infty]$.

## 1. Introduction

Several problems in extreme value theory and heavy tail analysis have as their core idea standard multivariate regular variation on a cone. Somewhat different theories and applications emerge by choosing different cones, but the common thread is that all of the problems have a reduction to regular variation of multivariate distributions on a particular cone. We consider the following three cones, labelled $\mathbb{C}$, and the associated three theories and application areas.

| Cone | Application |
|---|---|
| $\mathbb{C} = [\mathbf{0},\infty] \setminus \{\mathbf{0}\}$ | multivariate extreme value theory |
| $\mathbb{C} = (\mathbf{0},\infty]$ | hidden regular variation, coefficient of tail dependence; |
| $\mathbb{C} = [0,\infty] \times (0,\infty]$ | Conditioned limit theorems when one component is extreme. |

Table 1. Three theories stemming from standard multivariate regular variation on three different cones.

### 1.1. *Some background.*

We assume familiarity with the univariate theory of regularly varying functions. Fine surveys are available in [1–8]. Vector notation and commonly used symbols are listed in an appendix at the end and can be consulted as needed. We briefly review the relationship of vague convergence and one dimensional regular

Sidney Resnick's research was partially supported by ARO Contract W911NF-07-1-0078 and NSA Grant H98230-06-1-0069 at Cornell University. Most of the synthesis and writing of this paper took place during May–June 2007 when Sidney Resnick was visiting the Lehrstuhl für Mathematische Statistik, Technical University München and grateful acknowledgement is made for hospitality and support.

variation of distribution tails. When a function $U : \mathbb{R}_+ \mapsto \mathbb{R}_+$ is regularly varying with index $\rho \in \mathbb{R}$, we write $U \in RV_\rho$.

**1.1.1. Vague convergence.** Suppose $(\mathbb{E}, \mathcal{E})$ is a nice space (that is, locally compact with countable base). We set

$$M_+(\mathbb{E}) = \{\text{all Radon measures on } \mathbb{E}\}.$$

So $\mu \in M_+(\mathbb{E})$, means that $\mu$ is a measure on $\mathcal{E}$ and $\mu(K) < \infty$ for all compact sets $K \in \mathcal{E}$. Denote

$$C_K^+(\mathbb{E}) := \{f : f : \mathbb{E} \to \mathbb{R}_+,\ f \text{ continuous with compact support}\}.$$

If $\mu_n,\ \mu \in M_+(\mathbb{E})$ then $\mu_n \stackrel{v}{\to} \mu$ if for any $f \in C_K^+(E)$,

$$\mu_n(f) := \int_{\mathbb{E}} f d\mu_n \to \int_{\mathbb{E}} f d\mu =: \mu(f).$$

See [5; 7; 9–11].

The following result is from [5] and [7, page 62].

**Theorem 1.1:** *Regular variation of distribution function tails in dimension $d = 1$. Suppose $X \geq 0$ is a random variable with distribution function $F$ and set $\bar{F}(x) = 1 - F(x) = P[X > x]$. The following are equivalent:*

*(i)* $\bar{F} \in RV_{-\alpha},\ \alpha > 0$.

*(ii) There exists $b_n \to \infty$ such that*

$$n\bar{F}(b_n x) \to x^{-\alpha}, \quad x > 0.$$

*(iii) There exists $b_n \to \infty$ such that*

$$\mu_n(\cdot) = nP[\frac{X}{b_n} \in \cdot\,] \stackrel{v}{\to} \nu_\alpha(\cdot), \tag{1}$$

*in $M_+(0,\infty]$, where $\nu_\alpha \in M_+(0,\infty]$ satisfies $\nu_\alpha(x,\infty] = x^{-\alpha}$, for $x > 0$.*

**Remark 1:**

(a) We can always set

$$b(t) = \Big(\frac{1}{1-F}\Big)^{\leftarrow}(t) = F^{\leftarrow}(1 - \frac{1}{t}), \tag{2}$$

which we call the *quantile function* and then define $b_n = b(n)$.

(b) Note the use of the unusual state space $(0,\infty]$, which is a *one point uncompactification* of the compact set $[0,\infty]$. See [7, page 170]. This topology makes neighborhoods of $\infty$ relatively compact. This is required because vague convergence only controls behavior of measures on relatively compact sets and $(x,\infty]$ is a natural set when dealing with right tail problems.

## 2. Multivariate regular variation

Suppose $\mathbb{C}$ is a cone in $\mathbb{R}^d$, $d \geq 1$ so that $\boldsymbol{x} \in \mathbb{C}$ implies $t\boldsymbol{x} \in \mathbb{C}$, $\forall t > 0$. We assume that $\mathbf{1} \in \mathbb{C}$. A function $h : \mathbb{C} \to \mathbb{R}_+$ is regularly varying at $\infty$ with limit function $\lambda(\boldsymbol{x}) > 0$, if for all $\boldsymbol{x} \in \mathbb{C}$,

$$\lim_{t\to\infty} \frac{h(t\boldsymbol{x})}{h(t\mathbf{1})} = \lambda(\boldsymbol{x}) > 0\,, \quad \boldsymbol{x} \in \mathbb{C}\,. \tag{3}$$

For more details, see [7; 8].

We now consider an immediate implication: Fix $\boldsymbol{x} \in \mathbb{C}$ and define

$$U(t) = h(t\boldsymbol{x})\,, \quad t > 0\,.$$

Then $U \in RV_\rho$ for some $\rho \in \mathbb{R}$, since

$$\lim_{t\to\infty} \frac{U(ts)}{U(t)} = \lim_{t\to\infty} \frac{h(ts\boldsymbol{x})}{h(t\boldsymbol{x})} = \lim_{t\to\infty} \frac{h(ts\boldsymbol{x})}{h(t\mathbf{1})}\frac{h(t\mathbf{1})}{h(t\boldsymbol{x})} = \frac{\lambda(s\boldsymbol{x})}{\lambda(\boldsymbol{x})} = s^{\rho(\boldsymbol{x})}\,.$$

The exponent does not depend on $\boldsymbol{x}$ (see [7, page 169]) and therefore the limit function $\lambda$ is homogeneous:

$$\lambda(s\boldsymbol{x}) = s^\rho \lambda(\boldsymbol{x}). \tag{4}$$

### 2.1. *The polar coordinate transform.*

Let $\|\cdot\|$ be a norm on $\mathbb{R}^d$ and denote the unit sphere of $\mathbb{R}^d$ in this norm by $\aleph = \{\boldsymbol{x} \in \mathbb{R}^d : \|\boldsymbol{x}\| = 1\}$. Define $T : \mathbb{R}^d \setminus \{\mathbf{0}\} \mapsto (0,\infty) \times \aleph$ by

$$T(\boldsymbol{x}) = \Big(\|\boldsymbol{x}\|, \frac{\boldsymbol{x}}{\|\boldsymbol{x}\|}\Big) =: (r, \boldsymbol{a})\,. \tag{5}$$

The inverse map $T^{\leftarrow} : (0,\infty) \times \aleph \mapsto \mathbb{R}^d \setminus \{\mathbf{0}\}$ yields $T^{\leftarrow}(r, \boldsymbol{a}) = r\boldsymbol{a}$. Note that we do not define $T$ for $\boldsymbol{x} = \mathbf{0}$ or if $\boldsymbol{x}$ has infinite components.

### 2.2. *The one (or more) point uncompactification.*

The *one point uncompactification* identifies the compact sets in a compact space modified by the removal of one point. It is based on the following simple result. (See [7, page 170].)

**Proposition 2.1:** *Let $(\mathbb{X}, \mathcal{T})$ be a topological space with open sets $\mathcal{T}$. From $\mathbb{X}$, remove $\mathbb{D}$ to get*

$$\mathbb{X}^\sharp = \mathbb{X} \setminus \mathbb{D} = \mathbb{X} \cap \mathbb{D}^C$$

*with relative topology*

$$\mathcal{T}^\sharp = \mathcal{T} \cap \mathbb{D}^C = \mathcal{T} \cap \mathbb{X}^\sharp\,. \tag{6}$$

*Then the compact sets of $\mathbb{X}^\sharp$ are*

$$\mathcal{K}(\mathbb{X}^\sharp) = \{K \in \mathcal{K}(X) : K \cap \mathbb{D} = \phi\}\,. \tag{7}$$

So the compact sets in the new, modified space are the original compact sets prior to modification, provided they miss the deleted set $\mathbb{D}$.

**Example 2.2** Here are some examples of modified spaces that have been useful.

(i) $(0,\infty] = [0,\infty] \setminus \{0\}$ has compacta which are closed in $[0,\infty]$ and do not contain 0.
(ii) $[0,\infty]^d \setminus \{\mathbf{0}\} = [\mathbf{0},\infty] \setminus \{\mathbf{0}\}$ has compacta which are compact in $[0,\infty]^d$ and do not contain $\{\mathbf{0}\}$. We informally refer to such sets as bounded away from $\mathbf{0}$.
(iii) $(0,\infty]^d = [0,\infty]^d \setminus \cup_{i=1}^d \{t\boldsymbol{e}_i, t \geq 0\}$, where

$$\boldsymbol{e}_i = (0,\ldots 0,1,0,\ldots 0).$$

The compacta of this compact set modified by removal of the axes are closed subsets of $[0,\infty]^d$ not touching the axes. We refer to such sets as being neighborhoods of $\infty$.
(iv) $[0,\infty] \times (0,\infty] = [0,\infty]^2 \setminus \{t(1,0), t > 0\}$. The compact sets are closed subsets of $[0,\infty]^2$ not touching the horizontal axis.
(v) $(0,\infty] \times \aleph_+$, where $\aleph_+ = \{\boldsymbol{x} \geq \mathbf{0} : \|\boldsymbol{x}\| = 1\}$.

**2.3.** ***Multivariate regular variation of tail probabilities.***

We characterize regular variation of tail probabilities on the cone $\mathbb{C} = [\mathbf{0},\infty] \setminus \{\mathbf{0}\}$. Suppose $\boldsymbol{Z} \geq \mathbf{0}$ is a random vector with distribution $F$ concentrating on $[\mathbf{0},\infty)$. The following are equivalent but the symbols $\nu, b(\cdot), b_n$ may differ slightly in each instance.

(i) There exists a Radon measure on $\mathbb{C}$ called $\nu$ such that

$$\lim_{t\to\infty} \frac{1-F(t\boldsymbol{x})}{1-F(t\mathbf{1})} = \lim_{t\to\infty} \frac{P[t^{-1}\boldsymbol{Z} \in [0,\boldsymbol{x}]^c]}{P\big[t^{-1}\boldsymbol{Z} \in [\mathbf{0},\mathbf{1}]^c\big]} = \nu([\mathbf{0},\boldsymbol{x}]^c), \tag{8}$$

for all $\boldsymbol{x} \in [\mathbf{0},\infty) \setminus \{\mathbf{0}\}$ which are continuity points of the limit $\nu([\mathbf{0},\cdot]^c)$.
(ii) There exist $b(t) \to \infty$ and a Radon measure $\nu$ on $\mathbb{C}$ such that

$$tP[\boldsymbol{Z}/b(t) \in \cdot\,] \xrightarrow{v} \nu(\cdot), \quad t \to \infty, \text{ in } M_+(\mathbb{C}). \tag{9}$$

(iii) There exist $b_n \to \infty$, and a Radon measure $\nu$ on $\mathbb{C}$ such that in $M_+(\mathbb{C})$,

$$nP[\boldsymbol{Z}/b_n \in \cdot\,] \xrightarrow{v} \nu(\cdot), \quad n \to \infty.$$

(iv) Let $\aleph_+ = \aleph \cap [0,\infty]^d$. Then there exist a probability measure $S(\cdot)$ on Borel subsets of $\aleph_+$ called the angular measure, and a function $b(t) \to \infty$ such that with

$$(R,\boldsymbol{\Theta}) = \Big(\|\boldsymbol{Z}\|, \frac{\boldsymbol{Z}}{\|\boldsymbol{Z}\|}\Big)$$

we have

$$tP[\Big(\frac{R}{b(t)},\boldsymbol{\Theta}\Big) \in \cdot\,] \xrightarrow{v} c\nu_\alpha \times S \tag{10}$$

in $M_+\Big((0,\infty] \times \aleph_+\Big)$ or, equivalently,

$$nP[\Big(\frac{R}{b_n},\boldsymbol{\Theta}\Big) \in \cdot\,] \xrightarrow{v} c\nu_\alpha \times S\,. \tag{11}$$

The standard proofs are in [7, page 173].

We will need several extensions of the definition of multivariate regular variation given in (8), (9), or (10). We should consider other cones besides $[\mathbf{0},\infty] \setminus \{\mathbf{0}\}$. We should consider what happens if the components of $\boldsymbol{Z}$ are not normalized by the same scaling function $b(t)$. And finally, we should assess methods for analyzing dependence structure and how to statistically estimate such dependence.

### 2.4. *Some remarks on cones.*

Let $\overline{\mathbb{C}}$ be a closed compact cone in $[-\infty,\infty]^d$ with $\mathbf{0}\in\overline{\mathbb{C}}$. Suppose

$$\mathbb{C} := \overline{\mathbb{C}} \setminus (\{\mathbf{0}\}\cup\mathbb{D})$$

is still a cone. Compact subsets of $\mathbb{C}$ are the compact sets of $\overline{\mathbb{C}}$ not touching $\mathbf{0}$ or $\mathbb{D}$. The distribution of $\boldsymbol{Z}$ is regularly varying on $\mathbb{C}$ if there exists a scaling function $b(t)\to\infty$, and a limit Radon measure $\nu\in M_+(\mathbb{C})$ such that

$$tP[\frac{\boldsymbol{Z}}{b(t)}\in\cdot\,]\to\nu(\cdot),\quad \text{in } M_+(\mathbb{C})\,.$$

**Example 2.3** Consider the following examples.

| | Cone | $\mathbb{D}$ |
|---|---|---|
| 1. | $\mathbb{C}=[-\infty,\infty]^d\setminus\{\mathbf{0}\}$ | $\mathbb{D}=\emptyset$ |
| 2. | $\mathbb{C}=(0,\infty]^2$ | $\mathbb{D}=(0,\infty]\times\{0\}\cup\{0\}\times(0,\infty]$ |
| 3. | $\mathbb{C}=[0,\infty]\times(0,\infty]$ | $\mathbb{D}=(0,\infty]\times\{0\}$ |

Table 2. Three cones which are used for various applications.

Note that "Radon" means something different for each example since the notion of compact is different for each example. For instance:

- In the first cone, $\nu\{\boldsymbol{x}:\|\boldsymbol{x}\|>1\}<\infty$.
- Relative to the second cone, $\nu\{\boldsymbol{x}>\mathbf{0}:\|\boldsymbol{x}\|>1\}$ is not necessarily finite.
- For the third cone in the table, $\nu\{(x,y):\|(x,y)\|>1\}$ is also not necessarily finite.

**Example 2.4** Consider $d=1$ and suppose $\mathbb{C}=[-\infty,\infty]\setminus\{0\}$. Suppose $Z$ is a random variable with distribution $F$ and suppose $F$ has a regularly varying tail. There is no mass at $\{\pm\infty\}$ and regular variation means

$$tP[Z/b(t)\in\cdot\,]\to\nu(\cdot)\quad \text{in } M_+\Big([-\infty,\infty]\setminus\{0\}\Big)\,. \tag{12}$$

This implies that

$$tP[Z/b(t)>x]\to\nu(x,\infty],\quad x>0\,,\ t\to\infty, \tag{13}$$

and from the sequential form of regular variation, this yields

$$\nu(x,\infty]=c_+x^{-\alpha},\quad \alpha>0\,,\ c_+\geq 0\,. \tag{14}$$

Also, in a similar manner,

$$tP[Z/b(t)\leq -x]\to\nu[-\infty,-x],\quad x>0, \tag{15}$$

which implies

$$\nu[-\infty,-x]=c_-x^{-\alpha}\,. \tag{16}$$

Why is the $\alpha$ the same for both tails? There is only an issue if $c_+ > 0$, $c_- > 0$ since otherwise $\nu$ only concentrates on a half line emanating from the origin. So suppose $c_+ > 0$, $c_- > 0$. Then we have as $t \to \infty$,

$$t\bar{F}\big(b(t)x\big) \to x^{-\alpha}, \quad x > 0,$$

or equivalently after taking reciprocals

$$\frac{1}{t\Big(1 - F\big(b(t)x\big)\Big)} = \frac{U\big(b(t)x\big)}{t} \to x^{\alpha},$$

where $U = 1/(1-F)$ and we may take $U^{\leftarrow}(t) = b(t)$. But $U \in RV_\alpha$ implies $b \in RV_{1/\alpha}$.
Likewise for the other tail, as $t \to \infty$,

$$tF\big(b(t)(-x)\big) \to \nu[-\infty, -x], \quad x > 0,$$

leads to

$$\nu(-\infty, -x] = c_- x^{-\beta}$$

for some $\beta > 0$ which would imply $b(t) \in RV_{1/\beta}$ so we have $\beta = \alpha$.
We conclude the limit measure is of the form

$$\nu(dx) = c_+ \alpha x^{-\alpha-1} dx 1_{(0,\infty]}(x) + c_- |x|^{-\alpha-1} dx 1_{[-\infty,0)}(x)\,, \tag{17}$$

and also that

$$(a)\ P[|X| > x] \in RV_{-\alpha},$$
$$(b)\ \lim_{t\to\infty} \frac{P[X > t]}{P[|X| > t]} = \lim_{t\to\infty} \frac{tP[X/b(t) > 1]}{tP[|X|/b(t) > 1]} = \frac{\nu(1,\infty]}{\nu\{x : |x| > 1\}} = \frac{c_+}{c_+ + c_-} =: p,$$

and

$$(c)\ \lim_{t\to\infty} \frac{P[X < -t]}{P[|X| > t]} = \frac{c_-}{c_+ + c_-} =: q$$

where $p + q = 1$.
If $0 < \alpha < 2$, these are the necessary and sufficient conditions for convergence of sums of iid random variables with common distribution $F$ to a stable law. See [3; 7].

**Example 2.5** Suppose $d = 2$ and let $\boldsymbol{Z} = (Z^{(1)}, Z^{(2)})$ where $Z^{(1)}, Z^{(2)}$ are iid, and

$$P[Z^{(i)} > x] = x^{-1}, \quad x \geq 1.$$

The vector $\boldsymbol{Z}$ has two different regular variation properties on two different cones.

(a) $\boldsymbol{Z}$ has a regularly varying distribution on $[0,\infty]^2 \setminus \{\boldsymbol{0}\}$. To see this, note that for $i = 1, 2$,

$$tP[Z^{(i)} > tx] = t(tx)^{-1} = x^{-1}, \quad x > 0, \text{ for all } t \text{ big},$$

so we have marginal regular variation. For $\boldsymbol{x} = (x^{(1)}, x^{(2)})$, we have

$$\begin{aligned}\lim_{t\to\infty} tP\{[\boldsymbol{Z} \le t\boldsymbol{x}]^c\} &= \lim_{t\to\infty} tP\{[t^{-1}\boldsymbol{Z} \in [\boldsymbol{0}, \boldsymbol{x}]^c\} \qquad (18)\\ &= \lim_{t\to\infty} tP[Z^{(2)} > tx^{(2)}] + \lim_{t\to\infty} tP[Z^{(1)} > tx^{(1)}] - \lim_{t\to\infty} tP[Z^1 > tx^{(1)}]P[Z^{(2)} > tx^{(2)}]\\ &= (x^{(2)})^{-1} + (x^{(1)})^{-1} - 0 = \nu([0,x]^c). \qquad (19)\end{aligned}$$

This implies that $\nu$ concentrates on the axes and has no mass on $(0,\infty]^2$.

(b) $\boldsymbol{Z}$ has a distribution which is regularly varying on $(0,\infty]^2$: To see this, observe that

$$tP[\frac{\boldsymbol{Z}}{\sqrt{t}} \in (\boldsymbol{x}, \infty]] = \sqrt{t}P[Z^{(1)} > \sqrt{t}x^{(1)}]\sqrt{t}P[Z^{(2)} > \sqrt{t}x^{(2)}] \to (x^{(1)})^{-1}(x^{(2)})^{-1}\,. \qquad (20)$$

So on $[0,\infty]^2 \setminus \{\boldsymbol{0}\}$, $\boldsymbol{Z}$ has a distribution which is regularly varying with limit measure $\nu$ given in (19). On $(0,\infty]^2$, $\boldsymbol{Z}$ has a distribution which is regularly varying with limit measure $\nu_*(\boldsymbol{x},\infty] = (x^{(1)}x^{(2)})^{-1}$ and $b_*(t) = \sqrt{t}$.

### 2.5. *Form of the limit measure $\nu$.*

Suppose that $\nu$ is a Radon measure on Borel subsets of the cone $\mathbb{C}$, that $b(t) \to \infty$ is a scaling function and $F$ is a probability distribution such that

$$tF(b(t)\cdot) \xrightarrow{v} \nu(\cdot), \qquad \text{in } M_+(\mathbb{C}).$$

Recall vague convergence was defined in Section 1.1.1. For the cones we consider, we have the following properties:

(i) The scaling function $b(\cdot)$ satisfies $b \in RV_{1/\alpha}$ for some $\alpha > 0$.

(ii) The limit measure $\nu(\cdot)$ has a homogeneity property

$$\nu(t\cdot) = t^{-\alpha}\nu(\cdot), \quad \text{on } \mathbb{C}.$$

(iii) The scaling property (ii) of the limit measure $\nu$ is equivalent to a product form when using polar coordinates. To see this, suppose $\Lambda \in \mathcal{B}(\aleph \cap \mathbb{C})$. Then

$$\begin{aligned}\nu\{\boldsymbol{x} : \|\boldsymbol{x}\| > r\,,\ \frac{\boldsymbol{x}}{\|\boldsymbol{x}\|} \in \Lambda\} &= \nu\{\boldsymbol{x} : \|r^{-1}\boldsymbol{x}\| > 1\,,\ \frac{r^{-1}\boldsymbol{x}}{\|r^{-1}\boldsymbol{x}\|} \in \Lambda\}\\ &= \nu\{r\boldsymbol{y} : \|\boldsymbol{y}\| > 1\,,\ \frac{\boldsymbol{y}}{\|\boldsymbol{y}\|} \in \Lambda\}\\ &= r^{-\alpha}\nu\{\boldsymbol{y} : \|\boldsymbol{y}\| > 1\,,\ \frac{\boldsymbol{y}}{\|\boldsymbol{y}\|} \in \Lambda\}\,. \qquad (21)\end{aligned}$$

Define a measure $S$ on $\mathcal{B}(\aleph \cap \mathbb{C})$ by

$$S(\Lambda) := \nu\{\boldsymbol{y} : \|\boldsymbol{y}\| > 1\,,\ \frac{\boldsymbol{y}}{\|\boldsymbol{y}\|} \in \Lambda\}$$

and we conclude

$$\nu\{\boldsymbol{x} : \|\boldsymbol{x}\| > r\,,\ \frac{\boldsymbol{x}}{\|\boldsymbol{x}\|} \in \Lambda\} = r^{-\alpha}S(\Lambda)\,. \qquad (22)$$

Call $S$ the *angular or spectral measure.*

**Remark 1:** Note the following features.

(i) If $\aleph \cap \mathbb{C}$ is compact, $S$ is finite in which case we assume it is a probability measure.
(ii) If $\mathbb{C} = [0,\infty]^d \setminus \{\mathbf{0}\}$, then $\aleph \cap \mathbb{C} = \{\boldsymbol{x} \geq \mathbf{0} : \|\boldsymbol{x}\| = 1\}$ is compact.
(iii) If $\mathbb{C} = (0,\infty]^2$, then $\aleph \cap \mathbb{C} = \{\boldsymbol{x} > \mathbf{0} : \|\boldsymbol{x}\| = 1\}$ is not compact and $S$ is not necessarily finite.
(iv) If $\mathbb{C} = [0,\infty] \times (0,\infty]$,

$$\aleph \cap \mathbb{C} = \{(x,y) : 0 \leq x \leq \infty,\ \infty \geq y > 0,\ \|(x,y)\| = 1\}$$

is not compact.

*Expressing $\nu$ in terms of the angular measure $S$.* Suppose $B \in \mathcal{B}(\mathbb{C})$. Then integrating by means of polar coordinates we obtain

$$\nu(B) = \iint_{\{(r,\boldsymbol{a}):r\cdot\boldsymbol{a}\in B,\ r>0,\ \|\boldsymbol{a}\|=1\}} \alpha r^{-\alpha-1} dr S(d\boldsymbol{a}).$$

For example if $\mathbb{C} = [0,\infty]^d \setminus \{\mathbf{0}\}$,

$$\nu\big([\mathbf{0},\boldsymbol{x}]^c\big) = \iint_{\{r\boldsymbol{a})\leq x\}^c} \alpha r^{-\alpha-1} dr S(d\boldsymbol{a}). \tag{23}$$

Now

$$\begin{aligned}\{r\boldsymbol{a} \leq x\}^c =& \{(r,\boldsymbol{a}) : r > \frac{x^{(i)}}{a^{(i)}}, \quad \text{for some } i\}\\ =& \big\{(r,\boldsymbol{a}) : r > \bigwedge_{i=1}^{d} \frac{x^{(i)}}{a^{(i)}}\big\}.\end{aligned}$$

So integrate first in (23) with respect to $r$ via Fubini's theorem and then

$$\begin{aligned}\nu\big([\mathbf{0},\boldsymbol{x}]^c\big) &= \int_{\aleph\cap\mathbb{C}} \big(\bigwedge_{i=1}^{d} \frac{x^{(i)}}{a^{(i)}}\big)^{-\alpha} S(d\boldsymbol{a})\\ &= \int_{\aleph\cap\mathbb{C}} \bigvee_{i=1}^{d} \big(\frac{a^{(i)}}{x^{(i)}}\big)^{\alpha} S(d\boldsymbol{a}).\end{aligned} \tag{24}$$

### 3. Extreme value theory

Regular variation on the cone $[\mathbf{0},\infty] \setminus \{\mathbf{0}\}$ leads to the classical theory of multivariate extreme value distributions as limits of componentwise maxima of an iid sample. See [12], [5, Chapter 5] and [8].

Suppose $\boldsymbol{X}_*$ has distribution $F_*$ which is regularly varying on the cone $[\mathbf{0},\infty] \setminus \{\mathbf{0}\}$. We suppose the regular variation is in *standard form*, which means we may set $b(t) = t$ and there exists a Radon measure $\nu_*$ such that

$$tP[\frac{\boldsymbol{X}_*}{t} \in \cdot\,] \xrightarrow{v} \nu_*(\cdot), \quad \text{in } M_+\big([\mathbf{0},\infty] \setminus \{\mathbf{0}\}\big)\,. \tag{25}$$

From this we get the multivariate extreme value distributions in the following way, using the theory of extended regularly varying functions [2; 6]. Suppose $a(\cdot) > 0$, $b(\cdot) \in \mathbb{R}$ are measurable functions on $(0,\infty)$ satisfying,

$$\frac{b(tx) - b(t)}{a(t)} \to \psi(x), \quad x > 0, \, t \to \infty, \tag{26}$$

where $\psi \neq 0$, and $\psi$ is not constant. Then $\psi$ must be of the form

$$\psi(x) = \begin{cases} k\left(\frac{x^\rho - 1}{\rho}\right), & \rho \in \mathbb{R}\,,\ \rho \neq 0\,,\ x > 0 \\ k \log x & \rho = 0\,,\ x > 0 \end{cases} \tag{27}$$

for some $k \neq 0$. Furthermore,

(a) If $\rho > 0$, then $b(\cdot) \in RV_\rho$, and $b(t) \sim \frac{1}{\rho} a(t)$, and (26) really says that

$$\frac{b(tx)}{b(t)} \to x^\rho.$$

(b) If $\rho = 0$, then (26) defines $b(\cdot)$ as a $\Pi$-varying function. ($\Pi$-varying functions are a subclass of the slowly varying functions; see [2; 4–6; 8].) Also, $a(\cdot)$ is slowly varying and $b(t)/a(t) \to \infty$.

(c) If $\rho < 0$, then $b(\infty) := \lim_{t\to\infty} b(t)$ exists finite and

$$b(\infty) - b(t) \in RV_\rho\,, \quad b(\infty) - b(t) \sim \frac{1}{|\rho|} a(t)\,. \tag{28}$$

The inverted form of extended regular variation: If the function $b$ is nondecreasing with inverse $b^{\leftarrow}$, we can invert (26). Supposing for convenience that $k = 1$, the inversion of (26) yields

$$\frac{b^{\leftarrow}(b(t) + ya(t))}{t} \to \psi^{\leftarrow}(y) = (1 + \rho y)^{1/\rho}\,, \quad (1 + \rho y > 0). \tag{29}$$

Now assume (25) and suppose ($b^{(i)}, a^{(i)}$, $i = 1 \ldots d$ satisfy (26) with limits $\psi_i$ and parameters $\rho_i$ and that each $b^{(i)}$ is non-decreasing. Then

$$\begin{aligned} tP\Big\{\Big[\frac{b^{(i)}\big(X_*^{(i)}\big) - b^{(i)}(t)}{a^{(i)}(t)} \leq x^{(i)}\,, i = 1 \ldots d\Big]^c\Big\} \\ &= tP\Big\{\Big[\frac{X_*^{(i)}}{t} \leq b^{(i)\leftarrow}\big(a^{(i)}(t)x^{(i)} + b^{(i)}(t)\big)/t\,,\ i = 1 \ldots d\Big]^c\Big\} \\ &\to \nu_*\big(\{\boldsymbol{y} : y^{(i)} \leq \psi_i^{\leftarrow}(x^{(i)})\,,\ i = 1 \ldots d\}^c\big)\,. \end{aligned} \tag{30}$$

Where does the relation (30) lead? Set $X^{(i)} = b^{(i)}(X^{(i)})$ and suppose $b^{(i)}$ non-decreasing. Then (30) can be rephrased as

$$tP\{[X^{(i)} \leq a^{(i)}(t)x^{(i)} + b^{(i)}(t), i = 1, \ldots, d]^c\} \to \nu_*(\{\boldsymbol{y} : y^{(i)} \leq \psi_i^{\leftarrow}(x^{(i)});\, i = 1, \ldots, d\}^c)$$

or, replacing $t$ with $n$

$$\begin{aligned} nP\{[X^{(i)} \leq a^{(i)}(n)x^{(i)} + b^{(i)}(n);\ i = 1, \ldots, d]^c\} \\ &\to \nu_*(\{\boldsymbol{y} : y^{(i)} \leq \psi_i^{\leftarrow}(x^{(i)});\, i = 1, \ldots, d\}^c). \end{aligned} \tag{31}$$

If a sequence $p_n \to 0$ satisfies $np_n \to \ell \in (0,\infty)$, then also $n\big(-\log(1-p_n)\big) \to \ell$, and therefore (31) implies

$$n(-\log P[X^{(i)} \leq a^{(i)}(n)x^{(i)} + b^{(i)}(n);\ i=1,\dots,d] \to \nu_*(\{\boldsymbol{y} : y^{(i)} \leq \psi_i^{\leftarrow}(x^{(i)});\ i=1,\dots,d\}^c)$$

or exponentiating

$$\begin{aligned}\big(P[X^{(i)} \leq a^{(i)}(n)x^{(i)} + &b^{(i)}(n), i=1,\dots,d]\big)^n \\ &\to \exp\{-\nu_*(\{\boldsymbol{y} : y^{(i)} \leq \psi_i^{\leftarrow}(x^{(i)}), i=1,\dots,d\}^c)\}. \end{aligned} \tag{32}$$

Since

$$P[\bigvee_{j=1}^{n} \boldsymbol{X}_j \leq \boldsymbol{x}] = \big(P[\boldsymbol{X}_1 \leq \boldsymbol{x}]\big)^n,$$

when $\boldsymbol{X}_1,\dots,\boldsymbol{X}_n$ are iid random vectors and "$\leq$" is interpreted componentwise, we conclude, as $n \to \infty$,

$$P[\bigvee_{j=1}^{n}\Big(\frac{X_j^{(i)} - b^{(i)}(n)}{a^{(i)}(n)}\Big) \leq x^{(i)};\ i=1,\dots,d] \to \exp\{-\nu_*(\{\boldsymbol{y} : y^{(i)} \leq \psi_i^{\leftarrow}(x^{(i)});\ i=1,\dots,d\}^c)\}. \tag{33}$$

The limit is the form of the general max-stable distributions which are approximations for distributions of componentwise maxima of iid random vectors whose common distribution is in a maximal domain of attraction. See [5; 12; 13].

### 3.1. *Normalizing constants* $a(\cdot), b(\cdot)$.

Where do the functions $a(\cdot)$, $b(\cdot)$ satisfying (26) come from? Suppose $d=1$ and let $X_1,\dots,X_n$ be iid random variables with common distribution function $F(x)$ satisfying a one dimensional version of (33); that is, suppose there exist functions $a(t) > 0$, $b(t) \in \mathbb{R}$ and a disribution $G$ such that

$$P[\bigvee_{j=1}^{n} \frac{X_j - b(n)}{a(n)} \leq x] \to G(x)\,. \tag{34}$$

Then we say $F$ is in the *maximal domain of attraction of* $G$ and write $F \in MDA(G)$. The limit relation (34) leads to

$$P[\bigvee_{j=1}^{n} X_i \leq a(n)x + b(n)] = F^n\big(a(n)x + b(n)\big) \to G(x)\,,$$

so taking logarithms we get

$$n(-\log F\big(a(n)x + b(n)\big) \to -\log G(x).$$

This leads to

$$t(-\log F\big(a(t)x + b(t)\big) \to -\log G(x) \tag{35}$$

and

$$t(1 - F\big(a(t)x + b(t)\big) \to -\log G(x)$$

and then taking reciprocals,

$$\frac{1}{t}\Big(\frac{1}{1-F\big(a(t)x+b(t)\big)}\Big) \to \frac{1}{-\log G(x)}\,. \tag{36}$$

Invert to get

$$\frac{\Big(\frac{1}{1-F}\Big)^{\leftarrow}(yt)-b(t)}{a(t)} \to \Big(\frac{1}{-\log G}\Big)^{\leftarrow}(\,y), \quad y>0. \tag{37}$$

So we may set $b(t)=\big(\frac{1}{1-F}\big)^{\leftarrow}(t)$ and $a(t)=b(te)-b(t)$.

### 3.2. *Some conclusions.*

Assume again that $d>1$. We summarize the procedure outlined in this section.

(a) Start with a random vector $\boldsymbol{X}_*$ satisfying the standard form of regular variation (25) on $[\mathbf{0},\infty]\setminus\{\mathbf{0}\}$.
(b) Assume $b^{(i)}:(0,\infty)\mapsto(0,\infty)$ is non-decreasing and assume there exist $a^{(i)}(\cdot)>0$ such that for $i=1,\dots,d$ each $b^{(i)}$, $a^{(i)}$ pair satisfies (26).
(c) Transform $\boldsymbol{X}_*$ to $\boldsymbol{X}=\big(X^{(i)},\dots,X^{(d)}\big)$ by

$$X^{(i)}=b^{(i)}\big(X_*^{(i)}\big); \quad i=1,\dots,d.$$

(d) Then $\boldsymbol{X}$ satisfies (30) and if $(\boldsymbol{X}_1,\dots,\boldsymbol{X}_n)$ are iid and $\boldsymbol{X}_1\overset{d}{=}\boldsymbol{X}$ we get (33).
(e) We can choose

$$b^{(i)}(t)=\Big(\frac{1}{1-F_{(i)}}\Big)^{\leftarrow}(t) \tag{38}$$

where $F_{(i)}$ is a distribution in the maximal domain of attraction of a one dimensional extreme value distribution.

## 4. Hidden regular variation

Hidden regular variation is regular variation corresponding to the cone $\{\boldsymbol{x}\geq\mathbf{0}: x^{(i)}\wedge x^{(j)}>0,\text{ for some } 1\leq i,j\leq d\}$. Before reviewing this, there are several transitional topics that should be understood.

### 4.1. *A construction of multivariate regularly varying distributions.*

In this section $\mathbb{C}=[\mathbf{0},\infty]\setminus\{\mathbf{0}\}$. Suppose the random variable $R$ satisfies $R\geq 1$ and

$$P[R>r]=r^{-\alpha}, \quad r>1,$$

and $\boldsymbol{\Theta}$ is a random variable concentrating on $\aleph_+=\aleph\cap\mathbb{R}_+^d$ with distribution $S(d\boldsymbol{\theta})$. We assume $R$ and $\boldsymbol{\Theta}$ are independent, written $R\perp\!\!\!\perp\boldsymbol{\Theta}$. Then for $r>0$, $\Lambda\in\mathcal{B}(\aleph_+)$, we have

$$\begin{aligned} tP[\frac{R}{t^{1/\alpha}}>r\,,\ \boldsymbol{\Theta}\in\Lambda] =&tP[\frac{R}{t^{1/\alpha}}>r]P[\boldsymbol{\Theta}\in\Lambda]\\ \to& r^{-\alpha}S(\Lambda), \quad \text{as } t\to\infty. \end{aligned} \tag{39}$$

Note the same conclusion would hold if we only assumed $P[R>r]\in RV_{-\alpha}$ and replaced $t^{1/\alpha}$ by

$$b(t):=\Big(\frac{1}{1-P[R>\cdot]}\Big)^{\leftarrow}(t)\,.$$

Let $T$ be the polar coordinate transform: $T\boldsymbol{x}=(\|\boldsymbol{x}\|,\boldsymbol{x}/\|\boldsymbol{x}\|)=(r,\boldsymbol{a})$ and define

$$\boldsymbol{Z}=T^{\leftarrow}(R,\boldsymbol{\Theta}). \tag{40}$$

Then $\boldsymbol{Z}$ has a regularly varying distribution

$$tP[\frac{\boldsymbol{Z}}{b(t)}\in\cdot\,]\to\nu(\cdot) \tag{41}$$

where

$$\nu\circ T^{\leftarrow}=\alpha r^{-\alpha-1}dr S(d\boldsymbol{\theta}).$$

### 4.2. *The découpage de Lévy.*

The following result of Lévy helps to understand "peaks over threshold." See [5] for more detail on the proof.

**Theorem 4.1:** *Let $\{X_n,n>1\}$ be iid random elements of a nice space $(\mathbb{E},\mathcal{E})$. Fix $B\in\mathcal{E}$ with $P[X_1\in B]>0$ and define $\tau_0^{\pm}=0$ and*

$$\tau_i^{+}=\inf\{j>\tau_{i-1}^{+}:X_j\in B\},\quad \tau_i^{-}=\inf\{j>\tau_{i-1}^{-}:X_j\in B^c\}\,.$$

*Let $K_n=\sup\{i:\tau_i^{+}\le n\}$ be the (renewal) counting function for visits to $B$. Then*

$$\{X_{\tau_j^+}\}\,,\;\{X_{\tau_j^-}\}\,,\;\{K_n\}$$

*are independent with*

$$P[X_{\tau_1^+}\in A]=P[X_1\in A|X_1\in B]\,,\;A\subset B$$

*and*

$$P[X_{\tau_1^-}\in A]=P[X_1\in A|X_1\in B^c]\,,\;A\subset B^c\,.$$

*Further $\{K_n\}$ is a renewal counting function, $E(K_n)=nP[X_1\in B]$ and both $\{X_{\tau_j^+}\}$, $\{X_{\tau_j^-}\}$ are iid.*

### 4.3. *Peaks over threshold and the POT method.*

Peaks over threshold is an important concept. A good way to understand it is by means of the découpage de Lévy.

Suppose $\{\boldsymbol{Z}_j,j\ge1\}$ are iid $\mathbb{R}_+^d$-valued random vectors and set $B$ in the découpage to be

$$B=\{\boldsymbol{x}\in\mathbb{R}_+^d\;:\;||\boldsymbol{x}||>t\}\,.$$

Then for $s > 1$, $\Lambda \in \mathcal{B}(\aleph_+)$ and $ts > 1$, $\{\boldsymbol{Z}_{\tau_j^+}\}$ are iid and

$$\begin{aligned} P[R_{\tau_1^+} > ts\,,\ \boldsymbol{\Theta}_{\tau_1^+} \in \Lambda] =& P[R_1 > ts\,,\ \boldsymbol{\Theta}_1 \in \Lambda | R_1 > t] \\ =& \frac{P[R_1 > ts\,,\ \boldsymbol{\Theta}_1 \in \Lambda]}{P[R_1 > t]} \approx s^{-\alpha} S(\Lambda), \end{aligned} \tag{42}$$

where "$\approx$" means this is only true exactly as $t \to \infty$. The "peaks over threshold" or POT philosophy says *pretend* there is actual equality for some $t$. Then the method suggests,

(a) Consider the subsample of exceedances $\{\boldsymbol{Z}_{\tau_j^+}\,,\ j \geq 1\}$ with disribution in polar coordinates:

$$P[T(\boldsymbol{Z}_{\tau_1^+}) \in dr d\boldsymbol{\theta}] = \alpha r^{-\alpha-1} dr S(d\boldsymbol{\theta}), \quad r \geq 1,\, \boldsymbol{\theta} \in \aleph_+.$$

(b) Think of exceedances as coming from the construction of Section 4.1.

This method is popular for inference purposes but there is no obvious way to quantity error when replacing "$\approx$" with "=". See [14] for a convincing and clear discussion of advantages; there is also further material in [7].

### 4.4. *Asymptotic independence.*

Suppose $\boldsymbol{Z}$ is an $\mathbb{R}_+^d$-valued random vector with a distribution which is regularly varying on $[\boldsymbol{0},\infty] \setminus \{\boldsymbol{0}\}$ so that there exists $b(t) \to \infty$ such that

$$tP[\boldsymbol{Z}/b(t) \in \cdot\,] \to \nu(\cdot), \quad \text{in } M_+\big([\boldsymbol{0},\infty] \setminus \{\boldsymbol{0}\}\big)\,.$$

Then $\boldsymbol{Z}$ (or its distribution) is asymptotically independent if

$$\nu\{\boldsymbol{x} \in [\boldsymbol{0},\infty] \setminus \{\boldsymbol{0}\}\ :\ x^{(i)} \wedge x^{(j)} > 0, \quad \text{for some } 1 \leq i < j \leq d\} = 0 \tag{43}$$

or equivalently for $\boldsymbol{x} \geq \boldsymbol{0}$, $\boldsymbol{x} \neq \boldsymbol{0}$,

$$\nu\Big([\boldsymbol{0},\boldsymbol{x}]^c\Big) = \sum_{i=1}^{d} c_i (x^{(i)})^{-\alpha}\,. \tag{44}$$

Therefore, $\nu$ concentrates on the lines through $\boldsymbol{0}$ of the form $\{t\boldsymbol{e}_i\,,\ t > 0\}$, $i = 1,\dots,d$ where $\boldsymbol{e}_i = (0,\dots,0,1,0,\dots,0)$ with "1" in the $i$-th spot.

Recall the example $\boldsymbol{Z} = (Z^{(i)},\dots,Z^{(d)})$ where $Z^{(i)}$, $1 \leq i \leq d$ are iid and $P[Z^{(i)} > x) = x^{-\alpha}$, $x \geq 1$ leads to asymptotic independence, which is not surprising since independence ought to imply asymptotic independence.

There are two main reasons for the name *asymptotic independence* although this modifier has lots of other meanings in probability and statistics so we make no claim to universal appropriateness.

(i) Suppose $d = 2$ and $\boldsymbol{Z}$ possesses asymptotic independence. Then

$$\begin{aligned} \lim_{t\to\infty} P[Z^{(2)} > t | Z^{(1)} > t] = & \lim_{n\to\infty} \frac{P[Z^{(1)} > b(n), Z^{(2)} > b^{(n)}]}{P[Z^{(1)} > b(n)]} \\ = & \lim_{n\to\infty} cnP[Z^{(1)} > b(n), Z^{(2)} > b(n)] \\ =& c\nu(\boldsymbol{1},\infty] = 0 \end{aligned}$$

since $\nu$ concentrates on $\{(0,x): x>0\} \cup \{(x,0): x>0\}$. So given one component is large, the other tends not to be large; there is negligible probability they are both large.

(ii) The regular variation condition for the vector $\boldsymbol{Z}$ is equivalent to supposing $\{\boldsymbol{Z}_j, j \geq 1\}$ iid with $\boldsymbol{Z}_j \stackrel{d}{=} \boldsymbol{Z}$, $j \geq 1$, and

$$P[\bigvee_{j=1}^{n} \boldsymbol{Z}_j/b(n) \leq \boldsymbol{x}] \to e^{-\nu([\boldsymbol{0},\boldsymbol{x}]^c)} = P[\boldsymbol{Y} \leq \boldsymbol{x}],$$

for a limiting extreme value random vector $\boldsymbol{Y}$ with Frechet marginals. Asymptotic independence is equivalent to $Y^{(1)}, \dots, Y^{(d)}$, the components of $\boldsymbol{Y}$, being independent random variables [5, Chapter 5].

The name asymptotic independence is used frequently for many different concepts, so beware. Also note that it can sometimes be a surprising concept.

**Example 4.2** Consider the following two surprising examples and a third which is totally unsurprising.

(i) Suppose $U$ is $U(0,1)$ distributed and define

$$\boldsymbol{Z} = \Big(\frac{1}{U}, \frac{1}{1-U}\Big).$$

Then $\boldsymbol{Z}$ possesses asymptotic independence as can be verified by direct calculation. This example was pointed out in [15].

(ii) Let $(N_1, N_2)$ be bivariate normal with normal distribution

$$N\Big( \begin{pmatrix} 0 \\ 0 \end{pmatrix} \begin{pmatrix} 1 & \rho \\ \rho & 1 \end{pmatrix} \Big).$$

So $\mathrm{Corr}(N_1, N_2) = \rho$. Provided $\rho < 1$,

$$\boldsymbol{Z} = \Big(\frac{1}{\Phi(N_1)}, \frac{1}{\Phi(N_2)}\Big)$$

is asymptotically independent. See [16] or [5, Chapter 5].

(iii) Suppose $\boldsymbol{Z} = (Z^{(1)}, Z^{(2)})$ where $Z^{(1)} \perp\!\!\!\perp Z^{(2)}$ and for $i = 1, 2$,

$$P[Z^{(i)} > x] = x^{-1}, \quad x > 1.$$

This was already discussed.

### 4.5. *Hidden regular variation; definition and properties.*

Suppose we have two cones

$$\mathbb{C} = [\boldsymbol{0}, \infty] \setminus \{\boldsymbol{0}\}, \quad \mathbb{C}^0 = \mathbb{C} \setminus \bigcup_{i=1}^{d} \{t\boldsymbol{e}_i\,,\ t > 0\}.$$

We may also express $\mathbb{C}^0$ as

$$\mathbb{C}^0 = \{\boldsymbol{x} \in \mathbb{C} :\ \text{for some } 1 \leq i < j \leq d\,,\ x^{(i)} \wedge x^{(j)} > 0\}.$$

Suppose $\boldsymbol{Z} \geq \boldsymbol{0}$ has distribution $F$. Then if $F$ is regularly varying on $\mathbb{C}$ we say $F$ possesses hidden regular variation if it is also regularly varying on $\mathbb{C}^0$. Thus there exists a Radon measure $\nu$ on $\mathbb{C}$, and also there exists a Radon measure $\nu^0$ on $\mathbb{C}^0$ such that for scaling functions $b(t) \to \infty$ and $b^0(t) \to \infty$ we have

$$tP[\boldsymbol{Z}/b(t) \in \cdot\,] \xrightarrow{v} \nu(\cdot) \qquad \text{in } M_+(\mathbb{C}) \tag{45}$$

$$tP[\boldsymbol{Z}/b^0(t) \in \cdot\,] \xrightarrow{v} \nu^0(\cdot) \qquad \text{in } M_+(\mathbb{C}^0) \tag{46}$$

*and*, to ensure this is not a trivial concept, we assume

$$\frac{b(t)}{b^0(t)} \to \infty\,. \tag{47}$$

**4.5.1. Consequences and properties.** Here are some properties resulting from the definition. For fuller treatment, see [7; 17–19]. Hidden regular variation elaborates ideas of [20; 21].

(a) There exists $\alpha > 0$ such that for $t > 0$,

$$\nu(t\cdot) = t^{-\alpha}\nu(\cdot) \quad \text{ on } \mathbb{C}.$$

(b) There exists $\alpha^0 \geq \alpha$ such that for $t > 0$,

$$\nu^0(t\cdot) = t^{-\alpha^0}\nu^0(\cdot) \quad \text{ on } \mathbb{C}^0.$$

(c) Consequently, we have

$$b \in RV_{1/\alpha}, \quad b^0 \in RV_{1/\alpha^0}.$$

(d) $\aleph^0 := \aleph \cap \mathbb{C}^0$ is not compact.

(e) Suppose, as usual, that $T(\boldsymbol{x}) = (\|\boldsymbol{x}\|, \frac{\boldsymbol{x}}{\|\boldsymbol{x}\|})$ is the polar coordinate transform. Then for $c > 0$, we have

$$\nu \circ T^{-1} = c\nu_\alpha \times S \tag{48}$$

where $\nu_\alpha(x, \infty] = x^{-\alpha}$, for $x > 0$ and $S$ is a probability measure on $\aleph_+ = \aleph \cap \mathbb{C}$. Also for $c^0 > 0$,

$$\nu^0 \circ T^{-1} = c^0\nu_{\alpha^0} \times S^0 \tag{49}$$

where $\nu_{\alpha^0}(x, \infty] = x^{-\alpha^0}$ for $x > 0$, and $S^0$ is Radon on $\aleph^0$ but not necessarily finite.

(f) Regions of the form $(\boldsymbol{x}, \infty]$ are relatively compact in $\mathbb{C}^0$ and

$$\begin{aligned} tP[\bigwedge_{i=1}^d \frac{Z^{(i)}}{b^0(t)} > x] =& tP\Big[\frac{\boldsymbol{Z}}{b^0(t)} \in (x\mathbf{1}, \infty]\Big] \\ \to& \nu^0(x\mathbf{1}, \infty] = x^{-\alpha^0}\nu^0(\mathbf{1}, \infty]\,. \end{aligned} \tag{50}$$

With $x = 1$, the limit is $\nu^0(\mathbf{1}, \infty]$.

(g) Hidden regular variation implies asymptotic independence since if $\boldsymbol{z} > \mathbf{0}$

$$\begin{aligned} tP[\frac{\boldsymbol{Z}}{b(t)} > \boldsymbol{z}] \leq& tP\Big[\frac{\boldsymbol{Z}}{b(t)} > \bigwedge_{i=1}^d z^{(i)}\mathbf{1}\Big] \\ =& tP\Big[\frac{\boldsymbol{Z}}{b^0(t)} > \frac{b(t)}{b^0(t)} \bigwedge_{i=1}^d z^{(i)}\mathbf{1}\Big] \to 0, \end{aligned} \tag{51}$$

since $b(t)/b^0(t) \to \infty$ and $tP[\boldsymbol{Z}/b^0(t) \in \cdot\,] \to \nu^0(\cdot)$. So in contrast to (50) we get with $b(t)$ as the scaling function, that as $t \to \infty$,

$$tP[\bigwedge_{i=1}^{d} \frac{Z^{(i)}}{b(t)} > x] \to 0\,. \tag{52}$$

For $d = 2$, (52) gives

$$0 = \lim_{t\to\infty} P[Z^{(2)} > t | Z^{(1)} > t]$$

and if you are a religious copularian, this is

$$\lambda := \lim_{u\to 1} P[F_1(Z^{(2)}) > u | F_1(Z^{(1)}) > u] = 0$$

where $F_1$ is the marginal distribution.

(h) Choices for the scaling function $b(t)$ appearing in (45) include

$$b(t) = \big(\frac{1}{1-F_1}\big)^{\leftarrow}(t)$$

where $F_1$ is the distribution of $Z^{(1)}$ or alternatively,

$$b(t) = \big(\frac{1}{1-P[\bigvee_{i=1}^{d} Z^{(i)} > \cdot\,]}\big)^{\leftarrow}(t)\,.$$

Choices for $b^0(t)$ include

$$b^0(t) = \big(\frac{1}{1-P[\bigwedge_{i=1}^{d} Z^{(i)} > \cdot\,]}\big)^{\leftarrow}(t)\,.$$

This results from the next property.

(i) Regular variation on $\mathbb{C}$ implies $\bigvee_{i=1}^{d} Z^{(i)}$ has a regularly varying distribution tail

$$P[\bigvee_{i=1}^{d} Z^{(i)} > x] \in RV_{-\alpha}$$

while regular variation on $\mathbb{C}^0$ implies $\wedge_{i=1}^{d} Z^{(i)}$ has a regularly varying distribution tail,

$$P[\bigwedge_{i=1}^{d} Z^{(i)} > x] \in RV_{-\alpha^0}\,.$$

More precisely, we have the following result. See [7; 19].

**Theorem 4.3:** *Suppose $b(\cdot) \in RV_{1/\alpha}$, $b^0(\cdot) \in RV_{1/\alpha^0}$, $0 < \alpha \le \alpha^0$ and $\frac{b(t)}{b^0(t)} \to \infty$. Suppose also*

$$tP[\frac{Z^{(i)}}{b(t)} > x] \to x^{-\alpha}\,, \quad x > 0,\ i = 1, \ldots, d\,. \tag{53}$$

*Then $\boldsymbol{Z}$ possesses hidden regular variation iff*

(a) *Max-linear combinations have regularly varying tail probabilities,*

$$tP[\bigvee_{i=1}^{d} \frac{s^{(i)}Z^{(i)}}{b(t)} > x] \to c(\boldsymbol{s})x^{-\alpha}\,, \quad x > 0\,, t \to \infty. \tag{54}$$

*for* $\boldsymbol{s} \in [\boldsymbol{0}, \infty) \setminus \{\boldsymbol{0}\}$ *and some* $c(\boldsymbol{s}) > 0$*;*

*and*

(b) *Min-linear combinations have regularly varying tail probabilities,*

$$tP[\bigwedge_{i=1}^{d} \frac{a^{(i)}Z^{(i)}}{b^0(t)} > x] \to d(\boldsymbol{a})x^{-\alpha^0}\,, \quad x > 0\, t \to \infty, \tag{55}$$

*for some* $d(\boldsymbol{a}) > 0$ *where* $\boldsymbol{a} \in (\mathbb{C}^0)^{-1} = (\boldsymbol{0}, \infty] \setminus \bigcup_{i=1}^{d}\{t\boldsymbol{e}_i^{-1}, t > 0\}$ *and*

$$\boldsymbol{e}_i^{-1} = (\infty, \dots, \infty, 1, \infty, \dots, \infty).$$

(j) The last item suggests a diagnostic for detection: Assume $\boldsymbol{Z}_1, \dots, \boldsymbol{Z}_n$ is a random sample and compute the componentwise minima,

$$(\bigwedge_{i=1}^{n} Z_i^{(1)}, \dots, \bigwedge_{i=1}^{n} Z_i^{(d)})\,.$$

Estimate the tail index for the minima. If the problem has been standardized so that $\alpha = 1$ we seek evidence that $\hat{\alpha}^0 > 1$ which would be consistent with hidden regular variation.

**4.5.2. Hidden regular variation and marginal distributions.** Regular variation on $\mathbb{C}^0$ plus correct behavior of the marginal distributions implies hidden regular variation.

**Proposition 4.4:** *Suppose* $\boldsymbol{Z} \geq \boldsymbol{0}$ *satisfies*

$$Z^{(1)} \stackrel{d}{=} Z^{(2)} \stackrel{d}{=} \dots \stackrel{d}{=} Z^{(d)}$$

*and* $\boldsymbol{Z}$ *is regularly varying on* $\mathbb{C}^0$ *so that there exist* $b^0(\cdot) \in RV_{1/\alpha^0}$ *and* $\nu^0 \in M_+(\mathbb{C}^0)$ *such that*

$$tP[\frac{Z}{b^0(t)} \in \cdot\,] \xrightarrow{v} \nu^0(\cdot), \qquad \text{in } M_+(\mathbb{C}^0).$$

*Suppose additionally that one dimensional marginals have regularly varying tails; that is, there exists* $b(\cdot) \in RV_{1/\alpha}$ *with* $\alpha \leq \alpha^0$ *and* $b(t)/b^0(t) \to \infty$*, and*

$$tP[\frac{Z^{(i)}}{b(t)} > x] \to x^{-\alpha}\,, \quad x > 0,\ i = 1, \dots, d.$$

*Then* $\boldsymbol{Z}$ *is regularly varying on* $\mathbb{C} = [\boldsymbol{0}, \infty] \setminus \{\boldsymbol{0}\}$*.*

**Proof:** : For $\boldsymbol{x} > \boldsymbol{0}$ we have by the inclusion-exclusion formula,

$$
\begin{aligned}
tP\Big([\frac{\boldsymbol{Z}}{b(t)} \leq \boldsymbol{x}]^c\Big) =& t\sum_{i=1}^{d} P[\frac{Z^{(i)}}{b(t)} > x^{(i)}] - t\sum_{1\leq i<j\leq d} P[\frac{Z^{(i)}}{b(t)} > x^{(i)}\,,\ \frac{Z^{(j)}}{b(t)} > x^{(j)}] \\
&+ t\sum_{1\leq i<j<k\leq d} P[\frac{Z^{(i)}}{b(t)} > x^{(i)}\,, \frac{Z^{(j)}}{b(t)} > x^{(j)}\,, \frac{Z^{(k)}}{b(t)} > x^{(k)}] - \dots, \\
&\to \sum_{i=1}^{d} (x^{(i)})^{-\alpha} + 0, \quad (t \to \infty).
\end{aligned}
$$

The "0" comes from

$$
tP[\frac{Z^{(1)}}{b(t)} > x^{(1)}\,, \frac{Z^{(2)}}{b(t)} > x^{(2)}] = tP[\frac{Z^{(1)}}{b^0(t)} > x^{(1)}\frac{b(t)}{b^0(t)}\,, \frac{Z^{(2)}}{b(^0t)} > x^{(2)}\frac{b(t)}{b^0(t)}] \to 0,
$$

since $b(t)/b^0(t) \to \infty$. □

**4.5.3. Examples of hidden regular variation.** We now consider some examples to emphasize how hidden regular variation can arise in practice.

**Example 4.5**

(i) Recall Example 2.5. Suppose $\boldsymbol{Z} = (Z^{(1)}, Z^{(2)})$ with $Z^{(1)} \perp\!\!\!\perp Z^{(2)}$ and $P[Z^{(i)} > x] = x^{-1}$ for $x > 1$ and $i = 1, 2$. Then

$$
\begin{aligned}
&tP[\frac{\boldsymbol{Z}}{t} \in \cdot\,] \stackrel{v}{\to} \nu(\cdot), && \text{in } M_+(\mathbb{C}), \\
&tP[\frac{\boldsymbol{Z}}{\sqrt{t}} \in \cdot\,] \stackrel{v}{\to} \nu^0(\cdot), && \text{in } M_+(\mathbb{C}^0)
\end{aligned}
$$

where

$$
\nu(\boldsymbol{x}, \infty] = 0, \qquad \text{for } \boldsymbol{x} > \boldsymbol{0}, \text{ and} \tag{56}
$$

$$
\nu^0(\boldsymbol{x}, \infty] = \frac{1}{x^{(1)}x^{(2)}}\,, \qquad \text{for } \boldsymbol{x} > \boldsymbol{0}. \tag{57}
$$

(ii) Suppose $B$, $\boldsymbol{Y}$, $\boldsymbol{U}$ are independent with
  a) $P[B = 0] = P[B = 1] = \frac{1}{2}$.
  b) $\boldsymbol{Y} = (Y^{(1)}, Y^{(2)})$ where $Y^{(1)} \perp\!\!\!\perp Y^{(2)}$ and $P[Y^{(i)} > x] = x^{-1}$, for $x > 1$. Set $\alpha = 1$ and $b(t) = t$ to describe the regular variation of the distribution of $\boldsymbol{Y}$.
  c) $\boldsymbol{U}$ is multivariate regularly varying on $\mathbb{C}$ with index $\alpha^0$, $1 < \alpha^0 < 2$. Then there exists $b^0 \in RV_{1/\alpha^0}$ and a Radon measure $\nu^0 \in M_+(\mathbb{C})$ so that

$$
tP[\boldsymbol{U}/b^0(t) \in \cdot\,] \stackrel{v}{\to} \nu^0(\cdot) \not\equiv 0
$$

in $M_+(\mathbb{C})$. Here we really do mean $M_+(\mathbb{C})$ and not $M_+(\mathbb{C}^0)$.

Next, define

$$
\boldsymbol{Z} = B\boldsymbol{Y} + (1 - B)\boldsymbol{U}.
$$

For $\boldsymbol{x} > \mathbf{0}$, we have,

$$tP[\boldsymbol{Z}/b^0(t) > \boldsymbol{x}] = \tfrac{t}{2}P[\boldsymbol{Y} > b^0(t)\boldsymbol{x}] + \tfrac{t}{2}P[\boldsymbol{U} > b^0(t)\boldsymbol{x}] = I + II.$$

For $I$ we have because of the independent Pareto random variables that as $t \to \infty$,

$$\frac{t}{2}\frac{1}{b^0(t)^2 x^{(1)} x^{(2)}} \to 0$$

since $\big(b^0(t)\big)^2 \in RV_{2/\alpha^0}$ and $t/\big(b^0(t)\big)^2 \to 0$.

For $II$ we have

$$\mathrm{II} \to \frac{1}{2}\nu^0(\boldsymbol{x}, \infty]\,.$$

Also marginal distributions are $RV_{-\alpha} = RV_{-1}$, since

$$P[Z^{(1)} > x] = \tfrac{1}{2}P[Y^{(1)} > x] + \tfrac{1}{2}P[U^{(1)} > x] \in RV_{-1}$$

since the first term on the right is $RV_{-1}$ and the second is $RV_{-\alpha^0}$ with $\alpha^0 > 1$. This is enough to imply hidden regular variation by Proposition 4.4.

(iii) Here is a way to construct a class of examples.

Let $\Theta$ be a random variable with distribution $G$ concentrating on $(1, \infty)$ . Suppose $B$ is an independent Bernoulli random variable,

$$P[B = 0] = \tfrac{1}{2} = P[B = 1]\,.$$

Suppose $R \perp\!\!\!\perp \{B, \Theta\}$ and $P[R > r] = r^{-1}$, $r > 1$. Define

$$\boldsymbol{\Theta} = B(\Theta, 1) + (1 - B)(1, \Theta), \quad \boldsymbol{Z} = R\boldsymbol{\Theta}\,.$$

So the distribution of $\boldsymbol{\Theta}$ concentrates on the lines emanating from $(1, 1)$ in the horizontal and vertical directions.

We now give some properties of this construction.

a) We have

$$Z^{(1)} = R\Theta^{(1)}, \quad Z^{(2)} = R\Theta^{(2)},$$

and

$$\Theta^{(1)} = B\Theta + (1 - B)1, \quad \Theta^{(2)} = B1 + (1 - B)\Theta.$$

b) For the marginal distributions, we have,

$$P[Z^{(i)} > x] = P[R\Theta^{(i)} > x] = \int_1^\infty P[\Theta^{(i)} > x/r]\frac{dr}{r^2}$$

and setting $s = \frac{x}{r}$ we obtain

$$
\begin{aligned}
=&\tfrac{1}{x}\int_0^x P[\Theta^{(i)} > s]ds = \frac{1}{2x}[\int_0^1 2ds + \int_1^x \bar{G}(s)ds] \\
=&\frac{1}{2x}[2 + \int_1^x \bar{G}(s)ds] \,. \qquad (58)
\end{aligned}
$$

c) The random vector $\boldsymbol{Z}$ is regularly varying on $(0,\infty]^2$ with $\alpha^0 = 1$. Furthermore, let $\nu^0$ be the limit measure. It has angular measure $S^0$ and $S^0$ is finite iff

$$
\int_0^\infty \overline{G}(s)ds < \infty.
$$

To see this, note for $\nu^0$ we have for $\boldsymbol{x} > \boldsymbol{0}$,

$$
\nu^0(\boldsymbol{x}, \infty] = \lim_{t\to\infty} tP[\frac{\boldsymbol{Z}}{t} > \boldsymbol{x}] \,.
$$

Now

$$
\begin{aligned}
tP[\boldsymbol{Z}/t > \boldsymbol{x}] =& tP[R\Theta^{(1)} > tx^{(1)}\,, R\Theta^{(2)} > tx^{(2)}] \\
=& tP[R\Theta^{(1)} > tx^{(1)}\,, R\Theta^{(2)} > tx^{(2)}, B = 1] \\
&\quad + tP[R\Theta^{(1)} > tx^{(1)}\,, R\Theta^{(2)} > tx^{(2)}, B = 0] \\
=& \frac{t}{2}P[R\Theta > tx^{(1)},\, R > tx^{(2)}] + \frac{t}{2}P[R > tx^{(1)},\, R\Theta > tx^{(2)}] \\
=& \frac{t}{2}\int_{r > tx^{(2)}} \bar{G}(\frac{tx^{(1)}}{r} \vee 1)\frac{dr}{r^2} + \frac{t}{2}\int_{r > tx^{(2)}} \bar{G}(\frac{tx^{(2)}}{r} \vee 1)\frac{dr}{r^2}.
\end{aligned}
$$

Set $s = t/r$ and we obtain

$$
\begin{aligned}
=& \frac{1}{2}\int_0^{1/x^{(2)}} \bar{G}(sx^{(1)} \vee 1)ds + \frac{1}{2}\int_0^{1/x^{(1)}} \bar{G}(sx^{(2)} \vee 1)ds \\
=& \frac{1}{2x^{(1)}}\int_0^{x^{(1)}/x^{(2)}} \bar{G}(u \vee 1)du + \frac{1}{2x^{(2)}}\int_0^{x^{(2)}/x^{(1)}} \bar{G}(u \vee 1)du \,.
\end{aligned}
$$

Let, for instance, $x^{(1)} \to 0$ and the limit is infinite iff $\int_0^\infty \bar{G}(u)du = \infty$.

d) If $\bar{G} \in RV_{-\alpha}$, $\alpha < 1$, then $\boldsymbol{Z}$ is regularly varying on $\mathbb{C} = [0,\infty]^2 \setminus \{\boldsymbol{0}\}$ with index $\alpha$. The reason for this is that from item (iii)b and Karamata's theorem, we have

$$
P[Z^{(i)} > x] \in RV_{-\alpha}
$$

and the multivariate regular variation follows from Proposition 4.4.

**Remark 1:** Note the following.

(i) For Example (i), $\nu^0$ is infinite on $\mathbb{C}^0 \cap \{\boldsymbol{x} : \|\boldsymbol{x}\| > 1\}$. Hence $S^0$ is infinite on $\aleph^0$. This follows from (57) by, say, setting $x^{(1)} = 1$ and letting $x^{(2)} \to 0$.

(ii) For Example (ii), $\nu^0$ is finite on $\mathbb{C}^0 \cap \{\boldsymbol{x} : \|\boldsymbol{x}\| > 1\}$. Hence $S^0$ is finite on $\aleph^0 = \aleph \cap \mathbb{C}^0$. This follows because we assumed $\boldsymbol{U}$ is regularly varying on $\mathbb{C}$ and not just on $\mathbb{C}^0$.

(iii) For Example (iii), $\nu^0$ may be finite or infinite on $\mathbb{C}^0 \cap \{\boldsymbol{x} : \|\boldsymbol{x}\| > 1\}$, depending on whether the first moment of $G$ exists or not.

**4.5.4. The coefficient of tail dependence $\eta$ of Ledford and Tawn.** Hidden regular variation elaborates some ideas of [20; 21]. They define $\eta$ by assuming that $\boldsymbol{Z} = (Z^{(1)}, Z^{(2)})$ is a two dimensional random vector with non-negative components such that $Z^{(1)} \stackrel{d}{=} Z^{(2)}$ and

$$P[Z^{(1)} \wedge Z^{(2)} > t] \sim L(t)(P[Z^{(1)} > t])^{1/\eta}, \tag{59}$$

where $0 < \eta \le 1$ and $L \in RV_0$.

If $\boldsymbol{Z}$ possesses hidden regular variation, then (59) is true: Suppose as usual $b(\cdot)$ satisfies

$$tP[Z^{(1)} > b(t)] \sim 1.$$

Then (59) becomes

$$P[Z^{(1)} \wedge Z^{(2)} > b(t)] \sim L(b(t))t^{-1/\eta} = L'(t)t^{-1/\eta}, \tag{60}$$

where $L'(t) = L(b(t)) \in RV_0$. We know from hidden regular variation (See Theorem 4.3) that

$$P[Z^{(1)} \wedge Z^{(2)} > t] \sim t^{-\alpha^0}\ell(t),$$

for some function $\ell(t) \in RV_0$. So

$$P[Z^{(1)} \wedge Z^{(2)} > b(t)] \sim b(t)^{-\alpha^0}\ell(b(t)).$$

Now for some slowly varying function $\ell_*(t)$, $b(t) \sim t^{1/\alpha}\ell_*(t)$ and $\ell \circ b \in RV_0$, so

$$P[Z^{(1)} \wedge Z^{(2)} > b(t)] \sim t^{-\alpha^0/\alpha}\ell_{**}(t)\,, \tag{61}$$

for a slowly varying function $\ell_{**} \in RV_0$. From comparing (60) and (61), we conclude that

$$\eta = \frac{\alpha}{\alpha^0}\,.$$

## 5. Conditioned limit theorems

Suppose $\mathbb{C}^{1/2} = [0,\infty] \times (0,\infty]$ and $\boldsymbol{Z}$ is regularly varying on $\mathbb{C}^{1/2}$ meaning

$$tP[\frac{\boldsymbol{Z}}{b(t)} \in \cdot\,] \to \nu(\cdot), \quad \text{in } M_+(\mathbb{C}^{1/2}). \tag{62}$$

This implies

$$tP[\frac{Z^{(1)}}{b(t)} \in [0,\infty], \frac{Z^{(2)}}{b(t)} > x] = tP[\frac{Z^{(2)}}{b(t)} > x] \to \nu\big([0,\infty] \times (x,\infty]\big) = cx^{-\alpha}\,.$$

Suppose, as a normalization, that $c = 1$, so that,

$$\nu\big([0,\infty] \times (1,\infty]\big) = 1\,. \tag{63}$$

Then from (62), for $\Lambda \in \mathcal{B}(\mathbb{R}_+)$,

$$P[\frac{Z^{(1)}}{b(t)} \in \Lambda \Big| \frac{Z^{(2)}}{b(t)} > 1] = \frac{P[\frac{Z^{(1)}}{b(t)} \in \Lambda,\, \frac{Z^{(2)}}{b(t)} > 1]}{P[Z^{(2)}/b(t) > 1]}$$

$$\sim tP[\frac{Z^{(1)}}{b(t)} \in \Lambda,\, \frac{Z^{(2)}}{b(t)} > 1] \to \nu\big(\Lambda \times (1, \infty]\big)\,. \tag{64}$$

So regular variation on $\mathbb{C}^{1/2}$ implies a conditioned limit theorem: Given that the second component is large, the conditional distribution of the first variable is approximately determined by the limit measure in the definition of regular variation.

### 5.1. *General formulation.*

The previous formulation was for regular variation where each component could be normalized by the same scaling function. Regular variation on $\mathbb{C}^{1/2}$ naturally led to a conditioned limit theorem. How can we generalize this to make a more flexible class of statistical models? Here is an introduction; more detail and applications are found in [22–24].

We assume we have a two dimensional random vector $\boldsymbol{Z} = (X, Y)$. Suppose there exist $a(t) > 0$, $\alpha(t) > 0$, $\beta(t) \in \mathbb{R}$, $b(t) \in \mathbb{R}$ such that

$$tP\Big[\Big(\frac{X - \beta(t)}{\alpha(t)}, \frac{Y - b(t)}{a(t)}\Big) \in \cdot\Big] \stackrel{\nu}{\to} \mu(\cdot), \tag{65}$$

in $M_+\big([-\infty, \infty] \times (-\infty, \infty]\big)$ or equivalently

$$tP\Big[\frac{X - \beta(t)}{\alpha(t)} \le x, \frac{Y - b(t)}{a(t)} > y\Big] \to \mu\big([-\infty, x] \times (y, \infty]\big), \quad -\infty \le x \le \infty,\ y > -\infty. \tag{66}$$

We require that the limit measure $\mu$ satisfy the following constraints:

(a) For each fixed $y$, the distribution function $u[-\infty, x] \times (y, \infty]$ is not a degenerate distribution function in $x$.
(b) For $x \in [-\infty, \infty]$ and $y > -\infty$, we have

$$\mu([-\infty, x] \times (y, \infty]) < \infty.$$

### 5.2. *Standardizing the $Y$-variable*

Observe that (66) is in $M_+\big([-\infty, \infty] \times (-\infty, \infty]\big)$ and so we may insert the compact set $[-\infty, \infty]$ for the $X$-variable, effectively marginalizing to $Y$. Therefore,

$$tP[\frac{Y - b(t)}{a(t)} > y] \to \mu\big([-\infty, \infty] \times (y, \infty]\big), \quad y \in \mathbb{R}, \tag{67}$$

which entails that $Y \in D(G_\gamma)$ for some $\gamma \in \mathbb{R}$; that is, the distribution of $Y$ is in a maximal domain of attraction [5; 8; 25]. Note (67) is the defining property of univariate extreme value distributions (as in (35) in Section 3.1) so (67) entails (essentially by inverting (67)) that $b(\cdot)$ can be chosen as a monotone function and for some $\gamma \in \mathbb{R}$, as $t \to \infty$,

$$\frac{b(tx) - b(t)}{a(t)} \to c(\frac{x^\gamma - 1}{\gamma}), \quad x > 0, c \neq 0. \tag{68}$$

This allows us to standardize the $Y$-variable in the usual way. See [5, Chapter 5], or [8; 12]. Therefore write $Y^* = b^{\leftarrow}(Y)$ and

$$tP\Big[\frac{X-\beta(t)}{\alpha(t)} \le x\,,\frac{Y^*}{t} > y\Big) = tP\Big[\frac{X-\beta(t)}{\alpha(t)} \le x\,,\ \frac{Y-b(t)}{a(t)} > \frac{b(ty)-b(t)}{a(t)}\Big]$$
$$\to \mu\big([-\infty,x]\times(c(\frac{y^\gamma-1}{\gamma}),\infty]\big) =: \mu_{*/2}\big([-\infty,x]\times(y,\infty]\big), \quad x\in\mathbb{R}, y>0.$$

The subscript $*/2$ on the limit reminds us that the limit is obtained from the original $\mu$ by standardizing only the $Y$ variable.

Our conclusion: We can always standardize the $Y$-variable by making the transformation $Y \mapsto Y^* = b^{\leftarrow}(Y)$; then $Y^*$ need only be normalized by $t$. The meaning of the phrase *standardization* is that normalization by $t$ is adequate.

### 5.3. *Standardizing the X-variable.*

What about the $X$-variable? Suppose we know that $\beta$ is non-decreasing and write

$$tP[\frac{\beta^{\leftarrow}(X)}{t} \le x\,,\frac{Y^*}{t} > y] = tP[\frac{X-\beta(t)}{\alpha(t)} \le \frac{\beta(tx)-\beta(t)}{\alpha(t)}\,,\ \frac{Y^*}{t} > y]\,. \tag{69}$$

Provided

$$\frac{\beta(tx)-\beta(t)}{\alpha(t)} \to \psi_2(x) \not\equiv 0, \quad x>0, \tag{70}$$

we get with $X^* = \beta^{\leftarrow}(X)$ that for $x>0$, $y>0$,

$$tP[\frac{X^*}{t} \le x\,,\ \frac{Y^*}{t} > y] \to \mu\Big([-\infty,\psi_2(x)]\times\big(c(\frac{y^\gamma-1}{\gamma}),\infty]\Big) =: \mu_*\big([0,x]\times(y,\infty]\big), \tag{71}$$

and standardization is achieved.

**Proposition 5.1:** *[22] Suppose*

$$tP[\frac{X-\beta(t)}{\alpha(t)} \le x\,,\ \frac{Y^*}{t} > y] \to \mu_{*/2}([-\infty,x]\times(y,\infty]), \quad x\in\mathbb{R}, y>0, \tag{72}$$

*is the assumption where the $Y$-variable is standardized . We assume*

*(a) For each $y>0$, $\mu_{*/2}([-\infty,x]\times(y,\infty])$ is not a degenerate distribution in $x$; that is, it has more than one point of increase for the $x$-variable.*
*(b) The distribution of $Y^*$ is in the domain of attraction of $G_1$, the Frechet extreme value distribution with $\gamma=1$. We write this as $P[Y^*\le t]\in D(G_1)$. (This actually follows from (72).)*
*(c) The measure $\mu_{*/2}$ satisfies*

$$\mu_{*/2}\big([-\infty,\infty]\times(1,\infty]\big) = 1. \tag{73}$$

*Then there exist functions $\psi_1(c)$, $\psi_2(c)$, $c>0$ such that*

$$\lim_{t\to\infty}\frac{\alpha(tc)}{\alpha(t)} = \psi_1(c) \tag{74}$$

*and*

$$\lim_{t\to\infty} \frac{\beta(tc)-\beta(t)}{\alpha(t)} = \psi_2(c)\,, \tag{75}$$

*and the convergences in [74], [75] are locally uniform.*

**Proof:** Pick $c>0$ and assume $(x,1)$ and $(x,c^{-1})$ are continuity points of $\mu_{*/2}([-\infty,x]\times(y,\infty])$. (There are only countably many points not continuity points.) Define

$$H(x)=\mu_{*/2}([-\infty,x]\times(y,\infty]),\quad x\in\mathbb{R}, \tag{76}$$

and because of (73), $H(x)$ is a probability distribution function. Similarly, we define the conditional probability distribution

$$H_t\big(\alpha(t)x+\beta(t)\big)=P[\frac{X-\beta(t)}{\alpha(t)}\le x|Y^*>t]\,. \tag{77}$$

Then on the one hand,

$$\begin{aligned} H_t\big(\alpha(t)x+\beta(t)\big) =& \frac{P[\frac{X-\beta(t)}{\alpha(t)}\le x, Y^*>t]}{P[Y^*>t]} \\ \sim& tP[\frac{X-\beta(t)}{\alpha(t)}\le x\,,\ Y^*>t]\to \mu_{*/2}\big([-\infty,x]\times(1,\infty]\big)\,. \end{aligned} \tag{78}$$

On the other hand

$$\begin{aligned} H_t\big(\alpha(tc)x+\beta(tc)\big) =& P[\frac{X-\beta(tc)}{\alpha(tc)}\le x|\frac{Y^*}{t}>1]\sim tP[\frac{X-\beta(tc)}{\alpha(tc)}\le x, \frac{Y^*}{t}>1] \\ =& \frac{1}{c}\Big(tcP[\frac{X-\beta(tc)}{\alpha(tc)}\le x, \frac{Y^*}{tc}>\frac{1}{c}]\Big) \\ \to& \frac{1}{c}\mu_{*/2}\big([-\infty,x]\times(c^{-1},\infty]\big)=:H^{(c)}(x)\,. \end{aligned} \tag{79}$$

Comparing (78) and (79) and applying the convergence to types theorem [3; 26] we conclude that (74), (75) must hold and tha $H$ and $H^{(c)}$ are related by

$$H^{(c)}(x)=H\big(\psi_1(c)x+\psi_2(c)\big)\,. \tag{80}$$

□

### 5.4. *The limit measure is a product.*

When $\psi_2\not\equiv 0$, this allows standarization. However, it can occur that $\psi_2\equiv 0$ as shown in the following example [22; 24].

**Example 5.2** Let $N_1$, $N_2$ be iid $N(0,1)$ and define

$$(X,Y)=(\sqrt{1-\rho^2}N_1+\rho N_2, N_2)$$

so that $(X,Y)$ has a bivariate normal distribution

$$N\Big(\begin{pmatrix}0\\0\end{pmatrix},\begin{pmatrix}1\ \rho\\ \rho\ 1\end{pmatrix}\Big).$$

Define

$$\Phi(x)=\int_{-\infty}^{x}\frac{e^{-t^2/2}}{\sqrt{2\pi}}dt$$

for the normal cdf and then we may set (see, for instance, [5; 8])

$$a(t)=\frac{1}{\sqrt{2\log t}}, \tag{81}$$

$$b(t)=\big(\frac{1}{1-\Phi}\big)^{\leftarrow}(t)=\sqrt{2\log t}-\frac{\frac{1}{2}(\log\log t+\log 4\pi)}{\sqrt{2\log t}}+o(a(t)). \tag{82}$$

This choice of $a(\cdot)$ and $b(\cdot)$ ensures

$$tP[\frac{N_1-b(t)}{a(t)}>x]\to e^{-x},\quad x\in\mathbb{R}. \tag{83}$$

Furthermore [22; 24] for $(x,y)\in\mathbb{R}^2$,

$$tP[X-\rho b(t)\le x,\frac{Y-b(t)}{a(t)}>y]\to\Phi\big(\frac{x}{\sqrt{1-\rho^2}}\big)e^{-y}=\mu\big([-\infty,x]\times(y,\infty]\big), \tag{84}$$

or after standardization

$$tP[X-\rho b(t)\le x,\frac{Y^*}{t)}>y]\to\Phi\big(\frac{x}{\sqrt{1-\rho^2}}\big)y^{-1}=\mu_{*/2}\big([-\infty,x]\times(y,\infty]\big), \tag{85}$$

for $x\in\mathbb{R}$, $y>0$. It requires some calculation to verify (84) or (85) but the calculations are elementary. From (85), we see that

$$\alpha(t)=1\,,\ \beta(t)=\rho b(t)\,.$$

This implies

$$\frac{\beta(tc)-\beta(t)}{\alpha(t)}=\frac{\rho\big(b(tc)-b(t)\big)}{1}\,.$$

Either from (83) or from the explicit forms in (81) and (82) we get that the function $b\in\Pi$ with auxiliary function $a(\cdot)$, meaning [4; 5; 8]

$$\frac{b(tc)-b(t)}{a(t)}\to\log c,\quad c>0.$$

Thus we obtain

$$\frac{\beta(tc)-\beta(t)}{\alpha(t)}=\rho\Big(\frac{b(tc)-b(t)}{a(t)}\Big)a(t)\sim\rho\log c\cdot a(t)\to 0\equiv\psi(c)\,.$$

Also $\alpha(t) \equiv 1$ implies $\psi(c) \equiv 1$.

The circumstance $(\psi_1(c), \psi_2(c)) \equiv (1, 0)$ is associated with the limit measure being a product measure. We state the result [22] next.

**Proposition 5.3:** *The limit measure $\mu_{*/2}$ is a product measure*

$$\mu_{*/2} = H(\cdot) \times \nu_1 \qquad \text{where } \nu_1(x, \infty] = x^{-1}, x > 0, \tag{86}$$

*or equivalently*

$$\mu_{*/2}\big([-\infty, x] \times (y, \infty]\big) = H(x)y^{-1}, \quad y > 1, \tag{87}$$

*iff the functions $\psi_1, \psi_2$ in (74), (75) satisfy*

$$\psi_1(c) \equiv 1\,, \quad \psi_2(c) \equiv 0\,. \tag{88}$$

**Proof:** If $\mu_{*/2}$ is a product then from (79)

$$H^{(c)} = H$$

since (79) implies $\mu_{*/2}([-\infty, x] \times (1, \infty]) = H^{(c)}(x) = H(x)$. Therefore, $\psi_1 \equiv 1$, $\psi_2 \equiv 0$.

If $\psi_1 \equiv 1$ and $\psi_2 \equiv 0$, then $H^{(c)} = H$ and from (79)

$$\tfrac{1}{c}\mu_{*/2}([-\infty, x] \times (c^{-1}, \infty]) = H^{(c)}(x) = H(x)$$

or

$$\mu_{*/2}([-\infty, x] \times (c^{-1}, \infty]) = cH(x)$$

or

$$\mu_{*/2}([-\infty, x] \times (y, \infty]) = y^{-1}H(x)\,.$$

So $\mu_{*/2}$ is a product measure. □

## 6. Appendix

### 6.1. *Vector notation.*

Vectors are denoted by bold letters, capitals for random vectors and lower case for non-random vectors. For example: $\boldsymbol{x} = (x^{(1)}, \dots, x^{(d)}) \in \mathbb{R}^d$. Operations between vectors should always be interpreted componentwise so that for two vectors $\boldsymbol{x}$ and $\boldsymbol{z}$

$$\boldsymbol{x} < \boldsymbol{z} \text{ means } x^{(i)} < z^{(i)},\ i = 1, \dots, d, \qquad \boldsymbol{x} \le \boldsymbol{z} \text{ means } x^{(i)} \le z^{(i)},\ i = 1, \dots, d,$$

$$\boldsymbol{x} = \boldsymbol{z} \text{ means } x^{(i)} = z^{(i)},\ i = 1, \dots, d, \qquad \boldsymbol{z}\boldsymbol{x} = (z^{(1)}x^{(1)}, \dots, z^{(d)}x^{(d)}),$$

$$\boldsymbol{x} \vee \boldsymbol{z} = (x^{(1)} \vee z^{(1)}, \dots, x^{(d)} \vee z^{(d)}), \qquad \frac{\boldsymbol{x}}{\boldsymbol{z}} = \Big(\frac{x^{(1)}}{z^{(1)}}, \dots, \frac{x^{(d)}}{z^{(d)}}\Big),$$

and so on. Also, if $\boldsymbol{\alpha} = (\alpha^{(1)}, \dots, \alpha^{(d)}) \ge \mathbf{0}$, we write for $\boldsymbol{x} \ge \mathbf{0}$

$$\boldsymbol{x}^{\boldsymbol{\alpha}} = \big((x^{(1)})^{\alpha^{(1)}}, \dots, (x^{(d)})^{\alpha^{(d)}}\big).$$

Further, we define

$$\mathbf{0} = (0,\dots,0), \quad \mathbf{1} = (1,\dots,1), \quad \boldsymbol{\infty} = (\infty,\dots,\infty), \quad \boldsymbol{e}_i = (0,\dots,1,\dots,0),$$

where in $\boldsymbol{e}_i$ the "1" occurs in the $i$th spot. For a real number $c$, write as usual $c\boldsymbol{x} = (cx^{(1)},\dots,cx^{(d)})$. We denote the rectangles (or the higher dimensional intervals) by

$$[\boldsymbol{a},\boldsymbol{b}] = \{\boldsymbol{x} \in \mathbb{R}^d : \boldsymbol{a} \le \boldsymbol{x} \le \boldsymbol{b}\}.$$

Higher dimensional rectangles with one or both endpoints open are defined analogously, for example,

$$(\boldsymbol{a},\boldsymbol{b}] = \{\boldsymbol{x} \in \mathbb{R}^d : \boldsymbol{a} < \boldsymbol{x} \le \boldsymbol{b}\}.$$

**6.2.** ***Symbol and concept list.***

Here is a glossary of miscellaneous symbols and terminology used throughout.

| | |
|---|---|
| $RV_\rho$ | The class of regularly varying functions on $[0,\infty)$ with index $\rho \in \mathbb{R}$. |
| $D(\mathbb{E},\mathbb{R}^d))$ | The Skorokhod space of $\mathbb{R}^d$-valued càdlàg functions on the metric space $\mathbb{E}$ equipped with the $J_1$-topology. |
| $\epsilon_x$ | The probability measure consisting of all mass at $x$. |
| $f^{\leftarrow}$ | The left continuous inverse of a non-decreasing function $f$ defined by $f^{\leftarrow}(x) = \inf\{y : f(y) \ge x\}$. |
| $Leb$ | Lebesgue measure. |
| $\mathbb{E}$ | Often $[0,\infty]^d \setminus \{\mathbf{0}\}$ or $[-\infty,\infty] \setminus \{-\infty\}$. |
| Radon | Adjective applied to a measure to indicate the measure is finite on relatively compact sets. |
| $\stackrel{v}{\to}$ | Vague convergence of measures. |
| $M_+(\mathbb{E})$ | The space of non-negative Radon measures on $\mathbb{E}$. |
| $M_p(\mathbb{E})$ | The space of Radon point measures on $E$. |
| $\nu_\alpha$ | A measure on $(0,\infty]$ given by $\nu_\alpha(x,\infty] = x^{-\alpha}$, $\alpha > 0$, $x > 0$. |
| $\Rightarrow$ | Convergence in distribution. |
| $b(t)$ | Usually the quantile function of a distribution function $F(x)$, defined by $b(t) = F^{\leftarrow}(1 - \frac{1}{t})$ but usage can vary somewhat by context. |